\documentclass[12pt]{amsart}
\usepackage{a4,amsmath,amssymb,amsthm,eucal}
\usepackage{graphicx}
\usepackage{epstopdf}
\def\N{\mathbb{N}}
\def\R{\mathbb{R}}

\def\A{{\mathcal A}}

\def\K{{\mathcal K}}
\def\U{{\mathcal U}}

\def\cp{\mathrm{cap}\,}

\def\div{\mathrm{div}\,}
\def\dist{\mathrm{dist}\,}

\def\per{\mathrm{Per}}
\def\eps{\varepsilon}
\newtheorem{theorem}{Theorem}[section]
\newtheorem{definition}[theorem]{Definition}
\newtheorem{lemma}[theorem]{Lemma}

\theoremstyle{remark}
\newtheorem{remark}[theorem]{Remark}
\newtheorem{example}[theorem]{Example}

\numberwithin{equation}{section}

\author{G. Buttazzo\and A. Wagner}
\thanks{\it This research has been conceived during a visit of the first author to Department of Mathematics of RWTH Aachen University; he wishes to thank this institution for the warm and friendly atmosphere provided during all the visit.}
\title[Rescaled shape optimization problems]{On some rescaled shape optimization problems}
\begin{document}
\begin{abstract}
We consider Cheeger-like shape optimization problems of the form
$$\min\big\{|\Omega|^\alpha J(\Omega)\ :\ \Omega\subset D\big\}$$
where $D$ is a given bounded domain and $\alpha$ is above the natural scaling. We show the existence of a solution and analyze as $J(\Omega)$ the particular cases of the compliance functional $C(\Omega)$ and of the first eigenvalue $\lambda_1(\Omega)$ of the Dirichlet Laplacian. We prove that optimal sets are open and we obtain some necessary conditions of optimality.
\centerline{}
\end{abstract}
\maketitle
\section{Introduction}\label{sec1}
Many shape optimization problems are written in the form
\begin{equation}\label{shop}
\min\big\{F(\Omega)\ :\ |\Omega|=m,\ \Omega\subset D\big\}
\end{equation}
where $F$ is a suitable cost functional, $|\cdot|$ is the Lebesgue measure in $\R^N$, and $D$ represents a geometric constraint. For small values of $m$ often the optimal domains do not touch the boundary $\partial D$, which allows to obtain necessary conditions of optimality that lead in several cases to an explicit characterization of the solutions of \eqref{shop}.
A very well-known typical example is the isoperimetric problem
$$\min\big\{\per(\Omega)\ :\ |\Omega|=m,\ \Omega\subset D\big\}$$
where $\per(\cdot)$ is the De Giorgi perimeter; for small $m$ the solution is a ball, while in general for large $m$ a contact with the boundary $\partial D$ occurs, and the optimal domains have constant mean curvature in the free part. If the constraint on $|\Omega|$ changes, the perimeter scales as $|\Omega|^{1-1/N}$, in the sense that the quantity
\begin{equation}\label{perscale}
F(\Omega)=\frac{\per(\Omega)}{|\Omega|^\alpha}
\end{equation}
does not depend on $|\Omega|$ if $\alpha=1-1/N$. On the contrary, if $\alpha>1-1/N$ the minimum in the problem
$$\min\big\{F(\Omega)\ :\ \Omega\subset D\big\}$$
is reached on an optimal set $\Omega^*$ that touches $\partial D$. This is for instance the case of the Cheeger problem
\begin{equation}\label{cheeger}
\min\Big\{\frac{\per(\Omega)}{|\Omega|}\ :\ \Omega\subset D\Big\}
\end{equation}
where $\alpha=1$. It is known (see for instance \cite{ac07}, \cite{ccn07}) that for every bounded domain $D$ there exists an optimal Cheeger set $\Omega^*$ and this set is unique and convex whenever $D$ is convex. Moreover, in this case the boundary $\partial\Omega^*$ does not contain the points of $\partial D$ with too large mean curvature; more precisely, $\partial\Omega^*$ coincides with $\partial D$ if and only if
$$\|H\|_{L^\infty(\partial D)}\le\frac{\lambda(D)}{N-1}$$
where $H(x)$ is the mean curvature of $\partial D$ at $x$ and $\lambda(D)$ is the minimal value of the problem \eqref{cheeger}.

In the present paper we consider rescaled shape optimization problems for cost functionals related to elliptic equations, as functions of eigenvalues of the Dirichlet Laplacian or integral functionals depending on the solutions. More generally, we consider minimization problems of the form
\begin{equation}\label{rescaled}
\min\big\{M(\Omega)J(\Omega)\ :\ \Omega\subset D\big\}
\end{equation}
where the mappings $M$ and $J$ fulfill some rather general assumptions related to the variational $\gamma$-convergence on the family of quasi open sets. (we refer to \cite{bb} for a detailed presentation of this topic). In particular we do not require the monotonicity of $F(\Omega)=M(\Omega)J(\Omega)$. The powerful tools developed in the framework of this theory (see \cite{bb}, \cite{bdm91}, \cite{bdm93}, \cite{dm88}, \cite{dmm86}, \cite{dmm87}) allow us to obtain the existence of an optimal shape under quite general conditions.

Passing from the existence of an optimal domain $\Omega^*$ in the class of quasi open sets to the fact that $\Omega^*$ is actually an open set, and possibly to further regularity properties of $\Omega^*$, requires a very delicate analysis that is now available only for some particular problems of the form \eqref{rescaled}.

When the scaling factor of $F(\Omega)$ is above the scaling invariance, we have again that optimal sets $\Omega^*$ must touch the boundary $\partial D$. We analyze some particular cases in which this {\it``unnatural scaling''} allows to conclude that points of high mean curvature of $\partial D$ are not reached by $\Omega^*$.

We conclude the paper by a list of open questions which appear very natural.

\section{Preliminaries on capacity and related convergences}\label{sec2}
In the present paper we often use the notion of {\it capacity} of a subset $E$ of $\R^N$, defined by
$$\cp(E)=\inf\Big\{\int_{\R^N}|\nabla u|^2+u^2\,dx\ :\ u\in\U_E\Big\}\,,$$
where $\U_E$ is the set of all functions $u$ of the Sobolev space $H^1(\R^N)$ such that $u\ge1$ almost everywhere in a neighborhood of $E$. We resume here the main properties we shall use in the following; for all details about the capacity and the related convergences we will introduce subsequently we refer to the recent book \cite{bb}.

If a property $P(x)$ holds for all $x\in E$ except for the elements of a set $Z\subset E$ with $\cp(Z)=0$, we say that $P(x)$ holds {\it quasi-everywhere} (shortly {\it q.e.}) on $E$. The expression {\it almost everywhere} (shortly {\it a.e.}) refers, as usual, to the Lebesgue measure.

A subset $\Omega$ of $\R^N$ is said to be {\it quasi-open} if for every $\eps>0$ there exists an open subset $\Omega_\eps$ of $\R^N$, such that $\cp(\Omega_\eps{\scriptstyle\Delta}\Omega)<\eps$, where  ${\scriptstyle\Delta}$ denotes the symmetric difference of sets. Equivalently, a quasi-open set $\Omega$ can be seen as the set $\{u>0\}$ for some function $u$ belonging to the Sobolev space $H^1(\R^N)$. Notice that a Sobolev function is only defined quasi-everywhere, so a quasi-open set $\Omega$ does not change if we modify it by a set of capacity zero.

In this paper we fix a bounded open subset $D$ of $\R^N$ with a Lipschitz boundary and we consider the class $\A(D)$ of all quasi-open subsets of $D$. For every $\Omega\in\A(D)$ we denote by $H^1_0(\Omega)$ the space of all functions 
$u\in H^1_0(D)$ such that $u=0$ q.e.\ on $D\setminus\Omega$, endowed with the Hilbert space structure inherited from $H^1_0(D)$. In this way $H^1_0(\Omega)$ is a closed subspace of $H^1_0(D)$. If $\Omega$ is open, then the definition above of $H^1_0(\Omega)$ is equivalent to the usual one (see \cite{ah96}). If $\Omega\in\A(D)$ the linear operator $-\Delta$ on $H^1_0(\Omega)$ has a discrete spectrum and we denote by $\lambda_k(\Omega)$ the corresponding {\it eigenvalues}.

For every $\Omega\in\A(D)$ we consider the unique solution $w_\Omega\in H^1_0(\Omega)$ of the elliptic problem formally written as
\begin{equation}\label{eqw}
\left\{
\begin{array}{ll}
-\Delta w=1&\hbox{in }\Omega\\
w=0&\hbox{on }\partial\Omega
\end{array}\right.
\end{equation}
whose precise meaning has to be given through the weak formulation
$$\int_D\nabla w\nabla\phi\,dx=\int_D\phi\,dx\qquad\forall\phi\in H^1_0(\Omega).$$
The {\it compliance} functional $C(\Omega)$ is then defined as:
\begin{equation}\label{compli}
C(\Omega)=\int_\Omega w_\Omega\,dx.
\end{equation}
We introduce two useful convergences for sequences of quasi-open sets.

\begin{definition}\label{gamma}
A sequence $(\Omega_n)$ of quasi-open sets is said to $\gamma$-converge to a quasi-open set $\Omega$ if $w_{\Omega_n}\to w_\Omega$ in $L^2(\R^N)$.
\end{definition}

The following facts about $\gamma$-convergence are known (see \cite{bb}).

\begin{itemize}
\item The class $\A(D)$, endowed with the $\gamma$-convergence, is a metrizable and separable space, but it is not compact.

\item The $\gamma$-compactification of $\A(D)$ can be fully characterized as the class of all {\it capacitary measures} on $D$, that are Borel nonnegative measures, possibly $+\infty$ valued, that vanish on all sets of capacity zero.

\item The following maps are lower semicontinuous for the $\gamma$-convergence:
\begin{itemize}
\item for every integer $k$ the $k$-th eigenvalue $\lambda_k(\Omega)$ (they are actually $\gamma$-continuous);
\item the map $\Omega\mapsto\cp(D\setminus\Omega)$;
\item the compliance functional $C(\Omega)$ or more generally the integral functional
$${\ }\hskip1.5truecm C_{j,f}(\Omega)=\int_Dj(x,u_{\Omega,f},\nabla u_{\Omega,f})\,dx\quad\hbox{where}\quad
\left\{\begin{array}{ll}
-\Delta u=f&\hbox{in }\Omega\\
u=0&\hbox{on }\partial\Omega
\end{array}\right.$$
with $f\in H^{-1}(D)$ and $j(x,s,z)$ lower semicontinuous in $(s,z)$ and bounded from below by $-\alpha(x)-\beta(s^p+|z|^2)$ for suitable $\alpha\in L^1(\Omega)$ and $\beta\in\R$, where $p=2N/(N-2)$. The mapping above is actually $\gamma$-continuous if $j$ is a Carath\'eodory integrand with $|j(x,s,z)|\le\alpha(x)+\beta(s^p+|z|^2)$.
\end{itemize}
\end{itemize}
To overcome the lack of compactness of the $\gamma$-convergence, it is convenient to introduce another convergence, that we call $w\gamma$.

\begin{definition}\label{wgamma}
A sequence $(\Omega_n)$ of quasi-open sets is said to $w\gamma$-converge to a quasi-open set $\Omega$ if $w_{\Omega_n}\to w$ in $L^2(\R^N)$, and $\Omega=\{w>0\}$.
\end{definition}

We resume here the main facts about $w\gamma$-convergence (see \cite{bb}).

\begin{itemize}
\item The $w\gamma$-convergence is compact on the class $\A(D)$.

\item The $w\gamma$-convergence is weaker that the $\gamma$-convergence.

\item Every functional $F(\Omega)$ which is lower semicontinuous for the $\gamma$-convergence, and decreasing for the set inclusion, is lower semicontinuous for the $w\gamma$-convergence too. In particular, are $w\gamma$-lower semicontinuous:

\begin{itemize}
\item for every integer $k$, the map $\lambda_k(\Omega)$, and more generally the maps $\Phi(\lambda(\Omega))$ where $\lambda(\Omega)$ is the spectrum of the Dirichlet Laplacian in $\Omega$ and $\Phi:\R^\N\to[0,+\infty]$ is lower semicontinuous and nondecreasing (in each component);
\item the map $\cp(D\setminus\Omega)$;
\item the map $C_{j,f}(\Omega)$ when $f\ge0$ and $j(x,s,z)$ does not depend on $z$ and is decreasing in $s$.
\end{itemize}

\item The Lebesgue measure $|\Omega|$ is a mapping that is $w\gamma$-lower semicontinuous.
\end{itemize}

\section{Existence of optimal shapes}\label{sec3}
We consider cost functionals of the form $M(\Omega)J(\Omega)$ where $M$ and $J$ are defined on the class $\A(D)$ of all quasi open subsets of $D$, with values in $[0,+\infty]$. We assume:
\begin{eqnarray}
&&\hbox{$M$ and $J$ are nonnegative, and $J(D)>0$;}\label{posit}\\
&&\hbox{$J$ is $\gamma$-l.s.c. and nonincreasing with respect to the set inclusion;}\label{jsci}\\
&&\hbox{$M$ is $w\gamma$-l.s.c..\label{msci}}
\end{eqnarray}
Since $M$ or $J$ can take the value $+\infty$, in order to have the well posedness of the minimum problem we assume that
\begin{equation}\label{coerc}
\lim_{M(\Omega)\to0}M(\Omega)J(\Omega)=+\infty.
\end{equation}
With this assumption we may define the cost $M(\Omega)J(\Omega)=+\infty$ whenever $M(\Omega)=0$. The following existence result is now straightforward.

\begin{theorem}\label{exist}
Under the assumptions above the minimum problem
\begin{equation}\label{pbmin}
\min\big\{M(\Omega)J(\Omega)\ :\ \Omega\in\A(D)\big\}
\end{equation}
admits a solution $\Omega^*\in\A(D)$ and $M(\Omega^*)>0$.
\end{theorem}

\proof If $(\Omega_n)$ is a minimizing sequence, by the compactness of the $w\gamma$-convergence we may assume that $\Omega_n\to\Omega^*$ for some $\Omega^*\in\A(D)$, and $M(\Omega^*)>0$ by \eqref{coerc}. By the properties of $w\gamma$-convergence listed above the functional $J$ is $w\gamma$-l.s.c., as well as the product $M(\Omega)J(\Omega)$, which allows to conclude the proof.\endproof

The assumption \eqref{coerc} can be seen as a general condition which puts the minimum problem above the scaling invariance, as the following examples show.

\begin{example}\label{eigen}
If $M(\Omega)=|\Omega|^\alpha$ and $J(\Omega)=\lambda_1(\Omega)$, properties \eqref{posit}, \eqref{jsci}, \eqref{msci} are fulfilled. Since
$$|\Omega|^{2/N}\lambda_1(\Omega)\ge\lambda_1(B)$$
where $B$ is a ball of unit measure, condition \eqref{coerc} is fulfilled whenever $\alpha<2/N$. The same conclusion holds for $J(\Omega)=\lambda_k(\Omega)$ and more generally for $J(\Omega)=\Phi(\lambda(\Omega))$ where $\lambda(\Omega)$ is the spectrum of the Dirichlet Laplacian in $\Omega$, $\Phi:\R^\N\to[0,+\infty]$ is lower semicontinuous and nondecreasing (in each component), and $\Phi(\lambda)\ge c\lambda_1$ with $c>0$.
\end{example}

\begin{example}\label{compl}
If $M(\Omega)=|\Omega|^\alpha$ and $J(\Omega)=1/C(\Omega)$, properties \eqref{posit}, \eqref{jsci}, \eqref{msci} are fulfilled. Since
$$|\Omega|^{1+2/N}J(\Omega)\ge J(B)$$
where $B$ is a ball of unit measure, condition \eqref{coerc} is fulfilled whenever $\alpha<1+2/N$. The same conclusion holds for the functional
$$J(\Omega)=\int_\Omega j(x,u)\,dx$$
where $j(x,s)$ is a Carath\'eodory integrand, nondecreasing in $s$, satisfying $$|j(x,s)|\le\alpha(x)+\beta s^{2N/(N-2)}\qquad\hbox{and}\qquad j(x,s)\approx s^q\quad\hbox{as }s\to0.$$
In this case condition \eqref{coerc} is fulfilled whenever $\alpha<1+2q/N$.
\end{example}

\section{The case $M(\Omega)J(\Omega)=|\Omega|^\alpha/C(\Omega)$}\label{sec4}
In this section we will give a closer look to the case $F(\Omega):=M(\Omega)J(\Omega)$, where
\begin{equation}\label{functional}
M(\Omega)=|\Omega|^{\alpha},\qquad J(\Omega)=\big(C(\Omega)\big)^{-1},\qquad\alpha<1+\frac{2}{N},
\end{equation}
$C(\Omega)$ is the compliance functional defined in \eqref{compli} and $\Omega$ is assumed to vary in the class $\A(D)$ of all quasi-open subsets of $D$. Integrating by parts it is easy to see that
$$C(\Omega)=\min\Big\{\int_D|\nabla v|^2-2v\,dx\ :\ v\in H^1_0(\Omega)\Big\}.$$
The functional $C(\Omega)$ also allows a sup-formulation:
\begin{equation}\label{cvar}
C(\Omega)=\sup\big\{ R_C(v)\ :\ v\in H^1_0(\Omega)\big\}.
\end{equation}
where $R_C(v)$ denotes the quotient
$$R_C(v):=\frac{\Big(\int_\Omega v\,dx\Big)^2}{\int_\Omega|\nabla v|^2\,dx}.$$
Indeed a first variation of this quotient gives
$$\Delta v=-\big(R_C(v)\big)^{-1}\int_\Omega v\,dx$$
and since the quotient is invariant under scaling of $v$ we may assume 
\begin{equation}\label{Norm}
\int_\Omega v\,dx=R_C(v),\quad\hbox{i.e.}\quad\frac{\int_\Omega v\,dx}{\int_\Omega|\nabla v|^2\,dx}=1.
\end{equation}
This gives $\Delta v=-1$ and $C(\Omega)=\int_\Omega v\,dx$. Consequently we get
\begin{eqnarray}\label{cvar2}
\nonumber\big(C(\Omega)\big)^{-1}&=&\inf\big\{\big(R_C(v)\big)^{-1}\ :\ v\in H^1_0(\Omega)\big\}\\
&=&
\inf\Bigg\{\frac{\int_\Omega|\nabla v|^2\,dx}{\Big(\int_\Omega v\,dx\Big)^2}\ :\ v\in H^1_0(\Omega)\Bigg\}.
\end{eqnarray}
By the discussion in Example \ref{compl} we know that for $\alpha<1+\frac{2}{N}$ there exists an optimal domain $\Omega$ of $F(\Omega)$. In a first step we will show that $\Omega$ is open. The following remark is the key to that fact.

\begin{remark}\label{vonly}
We may reformulate the problem of finding the optimal set $\Omega\subset D$ as a problem in $v$ only. This is done by replacing the unknown set $\Omega$ by the positivity set of $v$.  Let
\begin{equation}\label{K(D)}
\K(D):=\{v\in H^1_0(D)\ :\ v\ge0\hbox{ a.e. in }D\}
\end{equation}
and for $v\in\K(D)$ define
$$F(v):=|\{v>0\}|^\alpha\,\frac{\int_{\{v>0\}}|\nabla v|^2\,dx}{\Big(\int_{\{v>0\}}v\,dx\Big)^2}.$$
Since $\int_{\{v>0\}}|\nabla v|^2\,dx=\int_D|\nabla v|^2\,dx$ and likewise the other integral, we get
\begin{equation}\label{jnew}
F(v)=|\{v>0\}|^{\alpha}\,\frac{\int_D|\nabla v|^2\,dx}{\Big(\int_D v\,dx\Big)^2}.
\end{equation}
Differently from \eqref{functional} this formulation involves only the unkonwn $v$ but not $\Omega$ any more. Clearly, any minimizer $v$ for \eqref{jnew} gives an domain $\Omega$ which is optimal in the sense of Theorem \ref{exist}: $\Omega=\{v>0\}$. Vice versa, any optimal domain $\Omega$ in the sense of Theorem \ref{exist} gives a minimizer $v$ for \eqref{jnew}.
\end{remark}

\begin{remark}\label{cnak}
Observe that 
$$\|v\|_{L^{2^*}(D)}\le c(N)\|\nabla v\|_{L^2(D)}$$
for any $v\in H^1_0(D)$, with $2^*=\frac{2N}{N-2}$. Thus
$$F(v)\ge c(N)|\{v>0\}|^\alpha\,\frac{\Big(\int_D|v|^{2^*}\,dx\Big)^{\frac{2}{2^*}}}{\Big(\int_Dv\,dx\Big)^2}$$
and by H\"older's inequality we obtain
\begin{equation}\label{star}
F(v)\ge c(N)|\{v>0\}|^{\alpha-1-\frac{2}{N}}.
\end{equation}
Since $\alpha<1+\frac{2}{N}$ we get $F(v)\to\infty$ for $|\{v>0\}|\to 0$, which gives property \eqref{coerc}. If $v$ is a minimizer of $F$ then in particular we have $F(v)\le K$ for some $K>0$. Thus \eqref{star} gives 
\begin{equation}\label{posmass}
|\{v>0\}|\ge c(N,\alpha,K)>0.
\end{equation}
If $\alpha\to 1+\frac{2}{N}$ it is easy to see that $c(N,\alpha,K)\to 0$.
\end{remark}

We compute now the first variation.
For $\delta>0$ we consider 
$$v_\delta(x):=\big(v(x)-\delta\varphi(x)\big)_+$$
with $\varphi\in C^\infty_0(D)$, $\varphi\ge0$. Clearly $v_\delta$ is in $\K(D)$ since the positive part of a function in $H^1_0(D)$ is still in $H^1_0(D)$. In particular, by minimality of $v$ we have $F(v)\le F(v_\delta)$ and with \eqref{Norm} this gives
$$F(v)\le|\{v>0\}|^\alpha R_C(v)^{-1}-2\delta\,\frac{|\{v>0\}|^\alpha}{\Big(\int_D v\,dx\Big)^2}\Big(
\int_D\nabla v\cdot\nabla\varphi\,dx-\int_D\varphi\,dx
\Big)+o(\delta).$$
Since $F(v)=|\{v>0\}|^\alpha R_C(v)^{-1}$ we get

\begin{theorem}\label{variation1} 
Let $v\in\K(D)$ be a minimum of $F(v)$. Then
\begin{equation}\label{var}
\Delta v + 1\ge0\quad\hbox{in }D
\end{equation}
in the distributional sense. Moreover, if $\{v>0\}$ contains an open subset $U$, then classical variation gives
$$\Delta v + 1= 0\quad\hbox{in }U.$$
\end{theorem}

We will prove the H\"older regularity for any minimizer $v$ of $F$.

\begin{theorem}\label{bound} 
Let $v$ be a minimizer of $F$. Then
\begin{equation}\label{bound1}
\| v\|_{L^{\infty}(D)}\le c(N)C(D)^{\frac{2}{N+2}}.
\end{equation}
\end{theorem}

\proof We use $(v-t)_{+}$ for some $t>0$ as an admissible test function in \eqref{var}. Thus we get
$$\int_D\nabla v\,\nabla (v-t)_+\,dx\le\int_D(v-t)_+\,dx.$$
Let $D(t)=\{v(x)>t\}$. Then we get
\begin{equation}\label{ineq1}
\int_{D(t)}|\nabla v|^2\,dx\le\int_{D(t)}(v-t)\,dx.
\end{equation}
H\"older inequality and Sobolev's imbedding give the inequalities
\begin{eqnarray*}
\int_{D(t)}(v-t)\,dx
&\le&
\Big(\int_{D(t)}|v-t|^{2^*}\,dx\Big)^{\frac{1}{2^*}}|D(t)|^{1-\frac{1}{2^*}}\\
&\le&
c(N)\Big(\int_{D(t)}|\nabla v|^2\,dx\Big)^{\frac{1}{2}}|D(t)|^{\frac{N+2}{2N}}.
\end{eqnarray*}
Thus we get
$$\int_{D(t)}|\nabla v|^2\,dx\ge
c(N)\Big(\int_{D(t)}(v-t)\,dx\Big)^2\|D(t)|^{-\frac{N+2}{N}}.$$
With this inequality we can estimate the left hand side of \eqref{ineq1}:
\begin{equation}\label{ineq2}
c(N)\Big(\int_{D(t)}(v-t)\,dx\Big)^2|D(t)|^{-\frac{N+2}{N}}
\le\int_{D(t)}(v-t)\,dx.
\end{equation}
Now recall the layer cake theorem (see e.g. Theorem 1.13 in \cite{ll01})
$$\int_{D(t)}(v-t)\,dx=\int_t^{+\infty}|D(s)|\,ds=:\tilde{D}(t).$$
With this \eqref{ineq2} reads as
$$c(N)\tilde{D}(t)\,|D(t)|^{-\frac{N+2}{N}}\le1.$$
Next we observe that
$$\tilde{D}'(t)=-|D(t)|.$$
Thus we finally get the differential inequality
$$c(N)\le-\big(\tilde{D}(t)\big)^{-\frac{N}{N+2}}\tilde{D}'(t)$$
which is equivalent to
$$c(N)\le-\frac{d}{dt}\big(\tilde{D}(t)^{\frac{2}{N+2}}\big).$$
Integration gives
$$c(N)\,t\le\tilde{D}(0)^{\frac{2}{N+2}}-\tilde{D}(t)^{\frac{2}{N+2}}$$
and hence
$$t\le c(N)\tilde{D}(0)^{\frac{2}{N+2}}
=c(N)\Big(\int_D v\,dx\Big)^{\frac{2}{N+2}}\le c(N)C(D).$$
as required.\endproof

Next we prove the H\"older continuity of $v$. This technique has been employed to similar problems in cases where the domain functional is monotone w.r.t. set inclusion (see \cite{wa05}, \cite{bawa07} ). The following lemma is crucial for what follows. A proof can be found e.g. in \cite{mazi97}.

\begin{lemma}\label{lDirichlet}
(Morrey's Dirichlet growth theorem). Let $v\in W^{1,p}(D)$, $1<p<N$. Suppose that there exist constants $0<c<+\infty$ and $\beta\in(0,1]$ such that for all balls $B_r(x_0)\subset D$
$$\int_{B_r(x_0)}|\nabla v|^p\,dx\le M\,r^{N-p+\beta p}.$$
Then $u\in C_{loc}^{0,\beta}(B)$. If the constant $M$ does not depend on $r$ and $\dist(x_0,\partial D)$ then $u\in C^{0,\beta}(\overline{D})$ and there exists a constant $c$ depending only on $D$, $N$ and $\beta$ such that
$$|v(x)-v(y)|\le c\,M|x-y|^\beta\qquad\forall x,y\in\overline{D}.$$
\end{lemma}

In order to apply the above lemma we shall also need

\begin{lemma}\label{lgrowth}
Let $\phi(t)$ be a nonnegative and nondecreasing function. Suppose that
$$\phi(r)\le\gamma\Big[\Big(\frac{r}{R}\Big)^\alpha+\delta\Big]\phi(R) +\kappa R^\beta$$
for all $0\le r\le R\le R_0$, where $\gamma$, $\kappa$, $\alpha$ and $\beta$ are positive constants with $\beta<\alpha$. Then there exist positive constants  
$\delta_0=\delta_0(\gamma,\alpha,\beta)$ and $c=C(\gamma,\alpha,\beta)$ such that if $\delta<\delta_0$, then
$$\phi(r)\le c\Big(\frac{r}{R}\Big)^\beta\big[\phi(R)+\kappa R^\beta\big]
\qquad\hbox{for all }0\le r\le R\le R_0.$$
\end{lemma}

For the proof of this Lemma we refer to \cite{gi83}, Lemma 2.1 in Chapter III. 

\begin{theorem}\label{holder1}
Let $v$ be a minimizer of $F$ such that $F(v)\le K$ for some $K>0$. Then $u\in C^{0,\beta}(\overline{D})$ and we have for any $0\le\beta<1$
\begin{equation}\label{bound2}
\|v\|_{C^{0,\beta}(\overline{D})}\le c(N,\alpha,K,D).
\end{equation}
\end{theorem}

We will first give the construction of an admissible comparison function and then prove some auxillary lemmas. We will use the notation $\Omega:=\{v>0\}$ and $\Omega_w:=\{w>0\}$. Let $x_0\in D$. Then there exists an $0<R<1$ such that $B_R(x_0)\subset D$. Consider the function
\begin{align}\label{constr.v}
w(x)=\begin{cases}
\hat{v}(x)&\text{if }x\in B_R(x_0)\\
v(x)&\text{if }x\in D\setminus B_R(x_0)
\end{cases}
\end{align}
where $\hat{v}$ is the solution of
\begin{equation}\label{hat}
\Delta\hat{v}+1=0\quad\hbox{in }B_R(x_0),\qquad
\hat{v}=v\quad\hbox{on }\partial B_R(x_0).
\end{equation}
By the strong maximum principle we have $\hat{v}>0$ in $B_R(x_0)$. Since $\Delta v+1\ge0$ in $D$ the maximum principle also gives 
\begin{equation}\label{comparison}
\hat{v}\ge v\quad\hbox{in }B_R(x_0).
\end{equation}
The function $w$ is admissible for variation, thus by minimality of $v$ we have $F(v)\le F(w)$. This is equivalent to
\begin{equation}\label{ineq3}
|\Omega|^\alpha R_C(v)^{-1}=|\Omega|^\alpha\,\frac{\int_D|\nabla v|^2\,dx}{\Big(\int_D v\,dx\Big)^2}
\le|\Omega_w|^\alpha\,\frac{\int_D|\nabla w|^2\,dx}{\Big(\int_D w\,dx\Big)^2}\;.
\end{equation}
We derive a local version of this inequality.

\begin{lemma}\label{local1}
Let $v$ be a minimizer of $F$ and $w$ defined as above. Then
\begin{eqnarray}\label{ineq4}
\int_{B_R(x_0)}|\nabla(v-\hat{v})|^2\,dx
&\le&
\frac{|\Omega\cup B_R(x_0)|^\alpha-|\Omega|^\alpha}{|\Omega\cup B_R(x_0)|^\alpha}\int_D|\nabla v|^2\,dx\\
\nonumber&&+2\int_{B_R(x_0)}(\hat{v}-v)\,dx.
\end{eqnarray}
\end{lemma}

\proof From \eqref{ineq3} and the definition of $w$ we deduce
\begin{eqnarray*}
&&|\Omega|^\alpha R_C(v)^{-1}\Big(\int_{D\setminus B_R(x_0)}v\,dx
+\int_{B_R(x_0)}\hat{v}\,dx\Big)^2\\
&&\qquad\qquad\le
|\Omega\cup B_R(x_0)|^\alpha
\Big(\int_{D\setminus B_R(x_0)}|\nabla v|^2\,dx+\int_{B_R(x_0)}|\nabla\hat{v}|^2\,dx\Big).
\end{eqnarray*}
Next we use $\hat{v}\ge v$ in $B_R(x_0)$ and the definition of $R_C(v)^{-1}$ to obtain
\begin{eqnarray*}
|\Omega|^\alpha\int_D|\nabla v|^2\,dx
&\le&
|\Omega\cup B_R(x_0)|^\alpha\int_D|\nabla v|^2\,dx\\
&&-|\Omega\cup B_R(x_0)|^\alpha
\int_{B_R(x_0)}|\nabla v|^2-|\nabla\hat{v}|^2\,dx.
\end{eqnarray*}
Rearranging terms then gives
\begin{equation}\label{ineq5}
\int_{B_R(x_0)}|\nabla v|^2-|\nabla\hat{v}|^2\,dx
\le
\frac{|\Omega\cup B_R(x_0)|^\alpha-|\Omega|^\alpha}{|\Omega\cup B_R(x_0)|^\alpha}\int_D|\nabla v|^2\,dx.
\end{equation}
Finally, since $\Delta\hat{v}=-1$ in $B_R(x_0)$ and $\hat{v}=v$ in $\partial B_R(x_0)$, we have
\begin{eqnarray*}
&&\int_{B_R(x_0)}|\nabla(v-\hat{v})|^2\,dx
=\int_{B_R(x_0)}|\nabla v|^2-2\nabla v\nabla\hat{v}+2|\nabla\hat{v}|^2-|\nabla\hat{v}|^2\,dx\\
&&\qquad=
\int_{B_R(x_0)}|\nabla v|^2-|\nabla\hat{v}|^2\,dx+2\int_{B_R(x_0)}\nabla\hat{v}\nabla(\hat{v}-v)\,dx\\
&&\qquad=
\int_{B_R(x_0)}|\nabla v|^2-|\nabla\hat{v}|^2\,dx+2\int_{B_R(x_0)}(\hat{v}-v)\,dx.
\end{eqnarray*}
Inserting this into \eqref{ineq5} gives
\begin{eqnarray*}
\int_{B_R(x_0)}|\nabla(v-\hat{v})|^2\,dx
&\le&
\frac{|\Omega\cup B_R(x_0)|^\alpha-|\Omega|^\alpha}{|\Omega\cup B_R(x_0)|^\alpha}\int_D|\nabla v|^2\,dx\\
&&+2\int_{B_R(x_0)}(\hat{v}-v)\,dx
\end{eqnarray*}
which shows the claim.\endproof

\begin{lemma}\label{mu}
Let $v$ be a minimizer of $F$ and $w$ defined as above. Then
\begin{equation}
\frac{|\Omega\cup B_R(x_0)|^\alpha-|\Omega|^\alpha}{|\Omega\cup B_R(x_0)|^\alpha}\le c(N,\alpha,K,D)R^N.
\end{equation}
\end{lemma}

\proof Set $X:=|\Omega|$ and $\delta:=|B_R(x_0)\setminus\Omega|$. Then
$$(X+\delta)^\alpha-X^\alpha=\alpha\int_0^\delta(X+t)^{\alpha-1}\,dt.$$
If $\alpha\le1$ we get 
$$\frac{(X+\delta)^\alpha-X^\alpha}{(X+\delta)^\alpha}\le\frac{\alpha}{X^\alpha}
\int_0^\delta(X+t)^{\alpha-1}\,dt\le\frac{\alpha\delta}{X}\;.$$
We use \eqref{posmass} to conclude that
$$\frac{(X+\delta)^\alpha-X^\alpha}{(X+\delta)^\alpha}\le\frac{\alpha\delta}{c(N,\alpha,K)}\;.$$
If $\alpha\ge1$ and w.l.o.g. $\delta\le X$ we get
$$\alpha\int_0^\delta(X+t)^{\alpha-1}\,dt\le\alpha2^{\alpha-1}\delta\,X^{\alpha-1}\le\alpha2^{\alpha-1}\delta\,|D|^{\alpha-1}.$$
Thus we obtain, using \eqref{posmass} again,
$$\frac{(X+\delta)^\alpha-X^\alpha}{(X+\delta)^\alpha}\le c(N,\alpha,K,D)R^N$$
as required.\endproof

In \eqref{ineq4} we can estimate the second term of the right hand side as follows. Using the fact that $\hat{v}=v$ in $\partial B_R(x_0)$ and Sobolev's imbedding we obtain
\begin{eqnarray*}
\int_{B_R(x_0)}(\hat{v}-v)\,dx
&\le&
c(N)\Big(\int_{B_R(x_0)}|\nabla(\hat{v}-v)|^2\,dx\Big)^{\frac{1}{2}}|B_R(x_0)|^{1-\frac{1}{2^*}}\\
&\le&
\frac{1}{4}\int_{B_R(x_0)}|\nabla(\hat{v}-v)|^2\,dx
+c(N) R^{N+2}
\end{eqnarray*}
Thus \eqref{ineq4} reads as 
$$\int_{B_R(x_0)}|\nabla(\hat{v}-v)|^2\,dx
\le c(N,\alpha,K,D)R^N\,\int_D|\nabla v|^2\,dx+c(N)R^{N+2}.$$
Next we observe that
$$\int_D|\nabla v|^2\,dx=-\int_D\Delta v\,v\,dx\le\int_D v\,dx=C(\Omega_v)\le C(D)$$
since $\Delta v+1\ge0$ (see \eqref{var}). Since $R<1$ we proved the following lemma.

\begin{lemma}\label{local2} 
Let $v$ be a minimizer of $F$ and $w$ defined as above. Then
\begin{equation}\label{ineq6}
\int_{B_R(x_0)}|\nabla(v-\hat{v})|^2\,dx\le c(N,\alpha,K,D) R^N.
\end{equation}
\end{lemma}

Finally we compare $\hat{v}$ with the harmonic function $h$ having the same boundary values on $\partial B_R(x_0)$.

\begin{lemma}\label{local3} 
Let $v$ be a minimizer of $F$ and $h$ defined as above. Then
\begin{equation}\label{ineq7}
\int_{B_R(x_0)}|\nabla(\hat{v}-h)|^2\,dx\le c(N) R^N.
\end{equation}
\end{lemma}

\proof From the weak formulation for $\hat{v}$ and $h$ and Sobolev's inequality we get
\begin{eqnarray*}
\int_{B_R(x_0)}|\nabla(\hat{v}-h)|^2\,dx
&=&
\int_{B_R(x_0)}|\hat{v}-h|\,dx\\
&\le&
c(N)\Big(\int_{B_R(x_0)}|\nabla(\hat{v}-h)|^2\,dx\Big)^{\frac{1}{2}}|B_R(x_0)|^{\frac{N+2}{2N}}.
\end{eqnarray*}
Rearranging terms gives
$$\int_{B_R(x_0)}|\nabla(\hat{v}-h)|^2\,dx\le c(N)R^{N+2}$$
which implies the claim for $0<R\le1$.\endproof

\noindent{\it Proof of Theorem \ref{holder1}.} For $0<r<R$ we have the well known growth estimate for harmonic functions
$$\int_{B_r(x_0)}|\nabla h|^2\,dx\le c(N)\Big(\frac{r}{R}\Big)^N\int_{B_R(x_0)}|\nabla h|^2\,dx.$$
Then we estimate
\begin{eqnarray*}
\int_{B_r(x_0)}|\nabla v|^2\,dx
&\le&
2\int_{B_r(x_0)}|\nabla(v-\hat{v}|^2\,dx
+2\int_{B_r(x_0)}|\nabla(\hat{v}-h)|^2\,dx\\
&&+
2\int_{B_r(x_0)}|\nabla h|^2\,dx\\
&\le&
2\int_{B_R(x_0)}|\nabla (v-\hat{v})|^2\,dx
+2\int_{B_R(x_0)}|\nabla (\hat{v}-h)|^2\,dx\\
&&+
2c(N)\Big(\frac{r}{R}\Big)^N\int_{B_R(x_0)}|\nabla h|^2\,dx.
\end{eqnarray*}
Applying Lemma \ref{local2} and Lemma \ref{local3} gives
$$\int_{B_r(x_0)}|\nabla v|^2\,dx
\le2c(N)\Big(\frac{r}{R}\Big)^N\int_{B_R(x_0)}|\nabla h|^2\,dx
+c(N,\alpha,C) R^N$$
We use now the fact that $h$ minimizes the Dirichlet integral among all functions in $H^1_0(B_R(x_0)$ having the same boundary values. Then
$$\int_{B_r(x_0)}|\nabla v|^2\,dx
\le2c(N)\Big(\frac{r}{R}\Big)^N\int_{B_R(x_0)}|\nabla v|^2\,dx
+c(N,\alpha,C) R^N.$$
Now we apply Lemma \ref{lgrowth}. This gives 
$$\int_{B_r(x_0)}|\nabla v|^2\,dx\le c(N)\Big(\frac{r}{R}\Big)^\beta
\int_{B_R(x_0)}|\nabla v|^2\,dx$$  
for all $0<\beta<N$. From Lemma \ref{lDirichlet} one has $u\in C^{0,\beta}(\overline{D})$ for all $0<\beta<1$.\endproof
\section{Some necessary conditions of optimality}\label{sec5}
In this section we find some necessary condition for the optimal domain $\Omega^*=\{v>0\}$, where $v\in\K(D)$ is a minimizer of $F(v)$ (see Section \ref{sec4} formula \eqref{jnew}).
Instead of first proving higher regularity for $v$ and $\partial\{v>0\}(=\partial\Omega^*)$ we assume $v\in C^1(\Omega^*)$ and $\partial\Omega^*\in C^{1,\beta}$ for some $0\le\beta<1$. From this we derive the desired optimality condition for the free boundary.

For any minimizer $v$ of \eqref{jnew} we consider a point $x_0\in\partial\Omega^* \cap D$. Let $R$ be such that $B_R(x_0)\subset D$. For any vector field $\eta\in C^\infty_0(B_R(x_0),\R^N)$ we define
$$v_\epsilon(x):= v(\tau_\epsilon^{-1}(x))\qquad\hbox{with}\qquad
\tau_\epsilon(x):=x+\epsilon\eta(x).$$
We expand with respect to $\epsilon$ and use the notation $D\eta=(\partial_i\eta_j)_{ij}$:
\begin{eqnarray*}
|\det D\tau_\epsilon|&=&1+\epsilon\,\div\eta+o(\epsilon)\\
(D\tau_\epsilon)^{-1}&=&Id-\epsilon\,D\eta+o(\epsilon).
\end{eqnarray*}
Set $\Omega_\epsilon=\{v_\epsilon>0\}$. We consider the functional
$$F(v_\epsilon)
=|\Omega_\epsilon|^\alpha\,\frac{\int_{\Omega_\epsilon}|\nabla v_\epsilon|^2\,dx}{\Big(\int_{\Omega_\epsilon}v_\epsilon\,dx\Big)^2}\;.$$
Then expansion above gives
\begin{eqnarray*}
&1)&|\Omega_\epsilon|^\alpha=|\Omega^*|^\alpha+\epsilon\alpha|\Omega^*|^{\alpha-1}\int_{\Omega^*}\div\,\eta\,dx+o(\epsilon)\\
&2)&\int_{\Omega_\epsilon}|\nabla v_\epsilon|^2\,dx=\int_{\Omega^*}|\nabla v|^2\,dx+\epsilon\Big(\int_{\Omega^*}-2\nabla v\cdot D\eta\cdot\nabla v+\div\,\eta\,|\nabla v|^2\,dx\Big)+o(\epsilon)\\
&3)&\Big(\int_{\Omega_\epsilon}v_\epsilon\,dx\Big)^{-2}=\Big(\int_{\Omega^*} v\,dx\Big)^{-2}
-2\epsilon\Big(\int_{\Omega^*}v\,dx\Big)^{-3}\int_{\Omega^*}\div\,\eta\,v\,dx+o(\epsilon),
\end{eqnarray*}
where $\frac{o(\epsilon)}{\epsilon}\to0$ as $\epsilon\to0$ and
$$\nabla v\cdot D\eta\cdot\nabla v=\sum_{i,j=1}^n\partial_iv\,\partial_i\eta_j\,\partial_jv.$$
For the integral in $2)$ we use partial integration. This gives
\begin{eqnarray*}
\int_{\Omega^*}-2\nabla v\cdot D\eta\cdot\nabla v+\div\,\eta\,|\nabla v|^2\,dx
&=&
2\int_{\Omega^*}\Delta v\,\eta\cdot\nabla v\,dx
-\int_{\partial\Omega^*}\eta\cdot\nu\,|\nabla v|^2\,dS\\
&=&
-2\int_{\Omega^*}\eta\cdot\nabla v\,dx
-\int_{\partial\Omega^*}\eta\cdot\nu\,|\nabla v|^2\,dS\\
&=&
2\int_{\Omega^*}\div\,\eta\,v\,dx
-\int_{\partial\Omega^*}\eta\cdot\nu\,|\nabla v|^2\,dS,
\end{eqnarray*}
since $\Delta v=-1$ in $\Omega^*$. From this we deduce
\begin{eqnarray*}
&&F(v_\epsilon)
=|\Omega_\epsilon|^\alpha\,\frac{\int_{\Omega_\epsilon}|\nabla v_\epsilon|^2\,dx}{\Big(\int_{\Omega_\epsilon}v_\epsilon\,dx\Big)^2}\\
&&=F(v)+\epsilon\Big(-2F(v)\Big(\int_{\Omega^*}v\,dx\Big)^{-1}\int_{\Omega^*}\div\,\eta\,v\,dx+\alpha F(v)|\Omega^*|^{-1}\int_{\Omega^*}\div\,\eta\,dx\Big)\\
&&\qquad+\epsilon\,|\Omega^*|^\alpha\Big(\int_{\Omega^*}v\,dx\Big)^{-2}\Big(2\int_{\Omega^*}\div\,\eta\,v\,dx-\int_{\partial\Omega^*}\eta\cdot\nu\,|\nabla v|^2\,dS\Big)+o(\epsilon).
\end{eqnarray*}
Again we use the fact that $\Delta v=-1$ in $\Omega^*$. Partial integration then leads to
$$F(v)=|\Omega^*|^\alpha\Big(\int_{\Omega^*}v\,dx\Big)^{-1}.$$
This simplifies the epression for the expansion.
$$F(v_\epsilon)=F(v)+
\epsilon F(v)\Big(\alpha|\Omega^*|^{-1}\int_{\Omega^*}\div\eta\,dx
-\Big(\int_{\Omega^*}v\,dx\Big)^{-1}\!\!\!\int_{\partial\Omega^*}\eta\cdot\nu\,|\nabla v|^2\,dS\Big)+o(\epsilon).$$
Since $C(\Omega^*)=\int_{\Omega^*}v\,dx$ we thus get
$$0\le\frac{F(v_\epsilon)-F(v)}{\epsilon}=\frac{F(v)}{C(\Omega^*)}\,
\int_{\partial\Omega^*}\eta\cdot\nu\,\Big(\alpha\frac{C(\Omega^*)}{|\Omega^*|}-|\nabla v|^2\Big)\,dS+\frac{o(\epsilon)}{\epsilon}.$$
Since $\eta\cdot\nu$ can have any sign we get $|\nabla v|^2=\alpha\frac{C(\Omega^*)}{|\Omega^*|}$ on $\partial\Omega^*\cap D$. If $x_0\in\partial\Omega^*\cap\partial D$ necessarily we have $\eta\cdot\nu\le0$. Thus we get $|\nabla v|^2\ge\alpha\frac{C(\Omega^*)}{|\Omega^*|}$ in $\partial\Omega^*\cap\partial D$. Thus we proved

\begin{theorem}\label{domainvar}
Let $v$ be a minimizer of $F$ and let $\Omega^*$ denote the set $\{v>0\}$. If $\partial\Omega^*\in C^{1,\beta}$ for some $0<\beta\le1$ and if $v\in C^1(\Omega^*)$, then necessarily the following conditions hold:
\begin{itemize}
\item $|\nabla v|^2=\alpha\frac{C(\Omega^*)}{|\Omega^*|}$ on $\partial\Omega^*\cap D$;
\item $|\nabla v|^2\ge\alpha\frac{C(\Omega^*)}{|\Omega^*|}$ on $\partial\Omega^*\cap\partial D$.
\end{itemize}
\end{theorem}

\begin{remark}One easily checks that the lower bound in Theorem \ref{domainvar} is strictly positive. Indeed, if we assume that $F(v)\le K$ for some $K>0$ we get
$$\frac{C(\Omega^*)}{|\Omega^*|}\ge\frac{|\Omega^*|^{\alpha-1}}{K}.$$
For $0\le\alpha\le1$ we get
$$\frac{C(\Omega^*)}{|\Omega^*|}\ge\frac{1}{K\,|D|^{1-\alpha}}.$$
For $1\le\alpha< 1+\frac{2}{N}$ inequality \eqref{posmass} implies
$$\frac{C(\Omega^*)}{|\Omega^*|}\ge\frac{c(N,\alpha,K)^{\alpha-1}}{K}.$$
\end{remark}

\begin{remark}If $\partial D$ has isolated conical points, but is smooth otherwise, it is well known, that the gradient of the solution $v$ of
$$-\Delta v=1\quad\hbox{in }\Omega,\qquad v=0\quad\hbox{in }\partial\Omega$$
vanishes in the conical points. Moreover there is a pointwise decay estimate for the gradient (see e.g. \cite{boko06} Theorem 3.11). Theorem \ref{domainvar} then shows, that the optimal domain cannot fill the entire set $D$.
\end{remark}

\section{Further remarks and problems}\label{sec6}
As we noticed in Example \ref{eigen} the existence Theorem \ref{exist} also applies to the case of cost functionals of the form $|\Omega|^\alpha\lambda_k(\Omega)$ where $\lambda_k(\Omega)$ is the $k$-th eigenvalue of the Dirichlet Laplacian in $\Omega$ and $\alpha<2/N$. Assuming that optimal domains are smooth enough, repeating computations similar to the one of Section \ref{sec5} we obtain the necessary conditions of optimality for the eigenfunction $u$
\begin{eqnarray*}
|\nabla u|^2&=&\alpha\frac{\lambda_k(\Omega^*)}{|\Omega^*|}\qquad\hbox{on}\quad\partial\Omega^*\cap D;\\
|\nabla u|^2&\ge&\alpha\frac{\lambda_k(\Omega^*)}{|\Omega^*|}\qquad\hbox{on}\quad\partial\Omega^*\cap\partial D.
\end{eqnarray*}
If $k=1$ it is possible to show (see \cite{bhp}) that optimal domains are actually open sets, whereas for $k\ge2$ this result, even is strongly expected, is not yet available.

Another class of problems occurs if we consider
$$M(\Omega)=\big(\per(\Omega)\big)^\alpha$$
where $\per(\Omega)$ is the perimeter of $\Omega$ in the sense of De Giorgi (see for instance \cite{afp}) and $\alpha$ is below the homogeneity threshold $2/(N-1)$. Even if the mapping $M(\Omega)$ is not in general $w\gamma$-l.s.c. it is possible to show (see \cite{bubuhe}) that the minimum problem
$$\min\big\{\big(\per(\Omega)\big)^\alpha\lambda_k(\Omega)\ :\ \Omega\subset D\big\}$$
admits a solution. The regularity of optimal domains and the corresponding necessary conditions of optimality have not yet been investigated.

We want to conclude the paper by pointing out some shape optimization problems for which the existence of a solution (though expected) is still unavailable. For a fixed $k\ge1$ we consider the optimization problem
$$\min\big\{C(\Omega)\lambda_k^\alpha(\Omega)\ :\ \Omega\in\A(D)\big\}$$
with $\alpha>1+N/2$, being the scaling invariance reached for $\alpha=1+N/2$. By 
the results of \cite{kj1} and \cite{kj2} we have
$$C(\Omega)\lambda_1^{1+N/2}(\Omega)\ge C(B)\lambda_1^{1+N/2}(B)$$
for any ball $B$, so that condition \eqref{coerc} is fulfilled whenever $\alpha>1+N/2$, by taking $M(\Omega)=C(\Omega)$ and $J(\Omega)=\lambda_k^\alpha(\Omega)$.

However, the $w\gamma$ l.s.c. condition \eqref{msci} fails for $C(\Omega)$, and so the existence Theorem \ref{exist} cannot be applied. It would be interesting to prove (or disprove) that an optimal domain $\Omega^*$ for the problem above exists.
\newline
\newline
{\bf{Acknowledgement:}} The authors like to thank Mark Ashbaugh for pointing out the work of M.T. Kohler-Jobin (\cite{kj1} and \cite{kj2}).

\bigskip
{\small
\begin{minipage}[t]{6.3cm}
Giuseppe Buttazzo\\
Dipartimento di Matematica\\
Universit\`a di Pisa\\
Largo B. Pontecorvo, 5\\
56127 Pisa - ITALY\\
{\tt buttazzo@dm.unipi.it}
\end{minipage}
\begin{minipage}[t]{6.3cm}
Alfred Wagner\\
Department of Mathematics\\
RWTH Aachen University\\
Templergraben 55\\
52062 Aachen - GERMANY\\
{\tt wagner@instmath.rwth-aachen.de}
\end{minipage}}
\end{document}